\UseRawInputEncoding
\documentclass[12pt,reqno]{amsart}
\usepackage{amsmath,amsfonts,amsthm,amsopn,amssymb}
\usepackage{cite,marginnote}
\usepackage{color,enumitem,graphicx}
\usepackage[colorlinks=true,urlcolor=blue,
citecolor=red,linkcolor=blue,linktocpage,pdfpagelabels,
bookmarksnumbered,bookmarksopen]{hyperref}
\usepackage[english]{babel}

\usepackage[left=2.9cm,right=2.9cm,top=2.8cm,bottom=2.8cm]{geometry}
\usepackage[hyperpageref]{backref}

\numberwithin{equation}{section}

\makeindex
\makeindex
\newtheorem{theorem}{Theorem}[section]
\newtheorem{definition}[theorem]{Definition}
\newtheorem{lemma}[theorem]{Lemma}
\newtheorem{corollary}[theorem]{Corollary}
\newtheorem{proposition}[theorem]{Proposition}
\newtheorem{remark}[theorem]{Remark}

\newcommand{\s}{\section}

\newcommand{\R}{\mathbb R}

\newcommand{\lab}{\label}
\newcommand{\bt}{\begin{theorem}}
\newcommand{\et}{\end{theorem}}
\newcommand{\bl}{\begin{lemma}}
\newcommand{\el}{\end{lemma}}
\newcommand{\bd}{\begin{definition}}
\newcommand{\ed}{\end{definition}}
\newcommand{\bc}{\begin{corollary}}
\newcommand{\ec}{\end{corollary}}
\newcommand{\bp}{\begin{proof}}
\newcommand{\ep}{\end{proof}}
\newcommand{\bx}{\begin{example}}
\newcommand{\ex}{\end{example}}
\newcommand{\bi}{\begin{exercise}}
\newcommand{\ei}{\end{exercise}}
\newcommand{\bo}{\begin{proposition}}
\newcommand{\eo}{\end{proposition}}
\newcommand{\br}{\begin{remark}}
\newcommand{\er}{\end{remark}}
\newcommand{\beq}{\begin{equation}}
\newcommand{\eeq}{\end{equation}}
\newcommand{\ba}{\begin{align}}
\newcommand{\ea}{\end{align}}
\newcommand{\bn}{\begin{enumerate}}
\newcommand{\en}{\end{enumerate}}
\newcommand{\bg}{\begin{align*}}
\newcommand{\bcs}{\begin{cases}}
\newcommand{\ecs}{\end{cases}}

\newcommand{\bean}{\begin{eqnarray*}}
\newcommand{\eean}{\end{eqnarray*}}

\def\N{\mathbb{N}}

\def\R{\mathbb{R}}

\def\bd{\mathrm{bd}\,}

\title[Normalized solution to mass super-critical problem]{Positive normalized solution to the Kirchhoff equation with general nonlinearities of mass super-critical}

\author[Q.~H.~He]{Qihan He}
\author[Z.~Y.~Lv]{Zongyan Lv}
\author[Y.~M.~Zhang]{Yimin Zhang}
\author[X.~X.~Zhong]{Xuexiu Zhong}

\address[Q.~H.~He]{\newline\indent College of Mathematics and Information Science
\newline\indent
Guangxi University,
\newline\indent
Nanning, 530003,PR China.}
\email{\href{mailto:heqihan277@163.com}{heqihan277@163.com}}

\address[Z.~Y.~Lv]{\newline\indent School of Mathematical Sciences
\newline\indent
Beijing Normal University
\newline\indent
Beijing 100875,PR China}
\email{\href{mailto:zongyanlv0535@163.com}{zongyanlv0535@163.com}}

\address[Y.~M.~Zhang]{\newline\indent Center for Mathematical Sciences
\newline\indent
Wuhan University of Technology
\newline\indent
Wuhan, 430070, PR China}
\email{\href{mailto:zhangym802@126.com}{zhangym802@126.com}}

\address[X.~X.~Zhong]{\newline\indent South China Research Center for Applied Mathematics and Interdisciplinary Studies
\newline\indent
South China Normal University
\newline\indent
 Guangzhou 510631, PR China}
\email{\href{mailto:zhongxuexiu1989@163.com}{zhongxuexiu1989@163.com}}

\thanks{The research was supported by the Natural Science Foundation of China  (No. 11801581, 1771127,12061012), Guangdong Basic and Applied Basic Research Foundation (2021A1515010034),Guangzhou Basic and Applied Basic Research Foundation(202102020225),the Fundamental Research Funds for the Central Universities (WUT: 2020IB011).}

\subjclass[2000]{}

\begin{document}

\begin{abstract}
In present paper, we study the normalized solutions $(\lambda_c, u_c)\in \R\times H^1(\R^N)$ to the following Kirchhoff problem
$$
-\left(a+b\int_{\R^N}|\nabla u|^2dx\right)\Delta u+\lambda u=g(u)~\hbox{in}~\R^N,\;1\leq N\leq 3
$$
satisfying the normalization constraint
$
\displaystyle\int_{\R^N}u^2=c,
$
which appears in free vibrations of elastic strings.  The parameters $a,b>0$ are prescribed as is the mass $c>0$. The nonlinearities $g(s)$ considered here are very general and of mass super-critical. Under some suitable assumptions, we can prove the existence of  ground state normalized solutions for any given $c>0$. After a detailed analysis via the blow up method, we also make clear the asymptotic behavior of these solutions as $c\rightarrow 0^+$ as well as $c\rightarrow+\infty$.
\end{abstract}

\maketitle

\s{Introduction}
\renewcommand{\theequation}{1.\arabic{equation}}

In this paper, we study the existence and the asymptotic behavior of the  ground state normalized solutions to the following Kirchhoff problem
\beq\lab{eq:Main-Problem}
-\left(a+b\int_{\R^N}|\nabla u|^2dx\right)\Delta u+\lambda u=g(u)~\hbox{in}~\R^N,\;1\leq N\leq 3
\eeq
satisfying the mass constraint
\beq\lab{eq:20210816-constraint}
\int_{\R^N}|u|^2 dx=c,
\eeq
where $a,b,c>0$ are prescribed.

The equation \eqref{eq:Main-Problem} was first proposed by Kirchhoff \cite{Ki} as an extension of the classical D'Alembert's wave equations for free vibration of elastic strings. It is worth mentioning that a functional analysis approach was introduced by J.L. Lions in the pioneer work \cite{JLLions1978}. After then, \eqref{eq:Main-Problem} has attracted extensive attentions.

In the equation \eqref{eq:Main-Problem}, if $\lambda$ is given, we call it {\it fixed frequency problem}. One can apply the variational method to look for the critical points of the associated energy functional $\Phi_\lambda: H^1(\R^N)\mapsto \R$, defined by
\beq\lab{eq:def-Phiu}
\Phi_\lambda[u]:=\frac{a}{2}\int_{\mathbb{R}^N}|\nabla u|^2dx+\frac{b}{4}\left(\int_{\mathbb{R}^N}|\nabla u|^2dx\right)^2+\frac{\lambda}{2}\int_{\R^N}|u|^2 dx-\int_{\mathbb{R}^N}G(u)dx£¬
\eeq
where $G(s):=\int_0^s g(t)dt$.
The existence of solutions is the first problem to be solved.
In \cite{G,LY}, positive ground state for equation (\ref{eq:Main-Problem}) is addressed. For  nontrivial solutions, we refer to \cite{ac,HL,Wu}.
The asymptotic behavior of the solutions is also a research interest of scholars.
In \cite{hz12,F,WT,wzz}, the authors considered \eqref{eq:Main-Problem} with very general nonlinearities, even involving the Sobolev critical exponent, and studied the existence, multiplicity and concentration behavior of positive solutions.  For the bounded states, we refer to \cite{H,HLP}. And the nodal solutions can be seen in \cite{DPS}.

We can also study problem \eqref{eq:Main-Problem} by searching for solutions satisfying the mass constraint \eqref{eq:20210816-constraint}. In such a point of view, the mass $c>0$ is prescribed, while the frequency $\lambda$ is unknown and will come out as a Lagrange multiplier. Hence, we call it {\it fixed mass problem} and the solution $u$ is called normalized solution. Normalized solutions of \eqref{eq:Main-Problem} can be searched as critical points of $I[u]$ constrained on $S_c$, where
 \beq\lab{eq:def-I}
I[u]:=\Phi_0[u]=\frac{a}{2}\int_{\mathbb{R}^N}|\nabla u|^2dx+\frac{b}{4}\left(\int_{\mathbb{R}^N}|\nabla u|^2dx\right)^2-\int_{\mathbb{R}^N}G(u)dx,
\eeq
and
\beq\lab{eq:def-Sc}
S_c:=\left\{u\in H^1(\R^N): \|u\|_2^2=c\right\}.
\eeq
As far as we know, for the study of normalized solution, the first work for equation \eqref{eq:Main-Problem} with $g(u)=|u|^{p-2}u$ is due to Ye. And it has been studied very well in a sequence of his work \cite{Y1,Y2,Y3}. Specifically, in \cite{Y1}, the author searched for minimizers to the following minimization problem:
\begin{equation}\label{eq1.2}
E_c:=\inf_{u\in S_c}I[u].
\end{equation}
By a scaling technique and applying the concentration-compactness principle, he proved that there exists $c_p^*\geq 0$, such that $E_c$ is attained if and only if $c>c_p^*$ with $0<p\leq 2+\frac{4}{N}$, or $c\geq c_p^*$  with $2+\frac{4}{N}<p<2+\frac{8}{N}$. The author also obtained that there is no minimizers for problem (\ref{eq1.2}) if $p\geq 2+\frac{8}{N}$. In particular, for the case of $2+\frac{8}{N}<p<2^*$, $E_c=-\infty$. However, the author could find a mountain pass critical point for the functional $I[u]$ constrained on $S_c$. In \cite{Y2}, Ye studied the case of $p=2+\frac{8}{N}$ and proved that there is a mountain pass critical point for the functional $I[u]$ on $S_c$ if $c>c^*$. Also, if $0<c<c^*$, the existence of minimizers for problem \eqref{eq1.2} was got by adding a new perturbation functional on the functional $I[u]$. Later on in \cite{Y3},  Ye also considered the asymptotic behavior of critical points of $I[u]$ on $S_c$ with $p=2+\frac{8}{N}$. Similar perturbation problems for $p=2+\frac{8}{N}$ can also be seen in \cite{zc}. By scaling technique and energy estimate, Zeng and Zhang \cite{zx} improved the results of \cite{Y1}. Precisely, they obtained the existence and uniqueness of minimizers of (\ref{eq1.2}) with $0<p<2+\frac{8}{N}$, and the existence and uniqueness of the mountain pass type critical points on the $L^2$ normalized manifold for $2+\frac{8}{N}<p<2^*$ or $p=2+\frac{8}{N}$ and $c>c^*$. We point out that some other existence results of normalized solutions of equation (\ref{eq:Main-Problem}) are  also considered in \cite{GW,lz} and the references therein.

According to these results, we know that the $L^2$-critical exponent $p=2+\frac{8}{N}$ is the threshold exponent related to the problem \eqref{eq:Main-Problem}. Precisely, the corresponding functional is bounded blow  if $p<2+\frac{8}{N}$, while there is no minimizers for problem \eqref{eq1.2} if $p\geq 2+\frac{8}{N}$.
However, there exists normalized solutions for equation (\ref{eq:Main-Problem}) if $p\geq 2+\frac{8}{N}$. Recently, Luo and Wang \cite{lw} studied the multiplicity of normalized solutions of equation (\ref{eq:Main-Problem}) in dimension $N=3$ with $g(u)=|u|^{p-2}u$ and $\frac{14}{3}<p<6$. The existence and blow-up analysis for equation (\ref{eq:Main-Problem}) with potential function can refer to \cite{gzz}. Indeed, for the homogeneous case that $g(s)=|s|^{p-2}s$,  by a scaling method, it is no hard to see the positive normalized solution  must be unique if exists.

It seems that no literatures involve the case of general nonlinearities, when it is nonhomogeneous and of mass super-critical.  It is also blank even for the existence. We note that in such a case, $E_c$ defined by \eqref{eq1.2} is not well defined. Hence, one can not find minimizers of $I[u]$ on $S_c$ any more. On the other hand, even though $g(s)=|u|^{p-2}u$ is homogenous, the asymptotic behavior of the solutions is also not studied before if $p>2+\frac{8}{N}$. Hence, our goal is to make some progresses in these respects.

In present paper, we shall introduce one more constraint, denoted by $\mathcal{P}_c$, see \eqref{eq:def-Pc}. We shall prove the new constraint $\mathcal{P}_c$ is natural, see Lemma \ref{lemma:20210807-xl2}. Then we can devote to search for the critical point of $I[u]$ on $\mathcal{P}_c$. We shall prove that the functional $I[u]$ possesses a mini-max structure, and the mini-max value denoted by $M_c$, coincides to the infimum of $I[u]$ constrained on $\mathcal{P}_c$, i.e.,
$$M_c:=\inf_{u\in S_c}\max_{t>0} I[t\star u]=\inf_{u\in \mathcal{P}_c}I[u],$$
where the fiber map $t\mapsto H^1(\R^N)$ is defined by
$\displaystyle(t\star u)(x):=t^{\frac{N}{2}}u(tx)$,
which preserves the $L^2$-norm.

\begin{definition}\lab{def-ground-state}
For any $c>0$, a solution $u\in H^1(\R^N)$ to \eqref{eq:Main-Problem} is  called  ground state normalized solution, or  least energy  normalized solution, if $u\in S_c$ and
$$I[u]=\min\left\{I[v]: v\in S_c~\hbox{and it solves \eqref{eq:Main-Problem} for some $\lambda\in \R$}\right\}.$$
\end{definition}
The nonlinearities considered in present paper are very general and of mass super-critical. For any given mass $c>0$, we shall prove that $M_c$ is attained by some $u_c\in H^1(\R^N)$, and corresponds some $\lambda_c\in \R$. One can check that a solution $u\in H^1(\R^N)$ to \eqref{eq:Main-Problem} satisfying $u\in S_c$ must belong to $\mathcal{P}_c$. Hence, $u_c$ is a ground state normalized solution to \eqref{eq:Main-Problem}. We also study the asymptotic behavior of $u_c$ as $c\rightarrow 0^+$ as well as $c\rightarrow+\infty$, via a blow up argument.

The paper is organized as follows. In the next section we state our main results.
Then in Section~\ref{sec:Preliminaries} we collect and prove a few basic facts about \eqref{eq:Main-Problem}. We shall prove the nimi-max structure of $I[u]$ and that $\mathcal{P}_c$ is a natural constraint. In Section \ref{Proof-th1}, we shall prove that $M_c$ is attained by some $u_c\in H_{rad}^{1}(\R^N)$, which is a positive decreasing function and the corresponding Lagrange multiplier $\lambda_c$ is also positive. Then Theorem \ref{thm:Main-Thoerem} will be proved there.
In Section \ref{sec:Properties-Mc}, we shall study some basic properties of $M_c$ and $\lambda_c$, which is crucial for our further study on the asymptotic behavior of $u_c$.
Theorems \ref{thm:20210814-th1} -\ref{thm:20210818-th2} will be proved in Section \ref{sec:aysmptotic-behavior}, via a blow up processes.

Throughout the paper we use the notation $\|u\|_p$ to denote  the $L^p$-norm. The notation $\rightharpoonup$ denotes weak convergence in $H^1(\R^N)$. Capital latter $C$ stands for positive constant which may depend on some parameters, whose precise value can change from line to line.

\section{Statement of results}
\renewcommand{\theequation}{2.\arabic{equation}}
To prove the existence of normalized solution to problem \eqref{eq:Main-Problem}, we devote to prove the achievement of $M_c$. And the major work is to construct a mini-max structure and the compactness analysis.
We firstly suppose that
\begin{itemize}
\item[(G1)]$g:\R\rightarrow \R$ is continuous and odd;
\item[(G2)] There exists some $(\alpha,\beta)\in \R_+^2$ satisfying
$$\begin{cases}
2+\frac{8}{N}<\alpha\leq \beta<2^*:=\frac{2N}{N-2}\quad&\hbox{if}~N=3,\\
2+\frac{8}{N}<\alpha\leq \beta<2^*:=+\infty \quad&\hbox{if}~N=1,2,
\end{cases}$$
such that
$$0<\alpha G(s)\leq g(s)s\leq \beta G(s)\;\hbox{for $s\neq 0$, where}~G(s)=\int_0^s g(t)dt.$$
\item[(G3)] The function defined by $\tilde{G}(s):=\frac{1}{2}g(s)s-G(s)$ is of class $C^1$ and
    $$\tilde{G}'(s)s\geq \alpha \tilde{G}(s), \forall s\in \R,$$
    where $\alpha$ is given by (G2).
\end{itemize}
Here comes our first main result:
\bt\lab{thm:Main-Thoerem}
Let $a,b>0$. Under the assumptions (G1)-(G3), for any given $c>0$, the Kirchhoff equation \eqref{eq:Main-Problem} has a  ground state normalized solution $(\lambda_c, u_c)$ with $\lambda_c>0$ and $u_c\in H_{rad}^{1}(\R^N)$.
\et

To study the asymptotic behavior of $u_c$, we suppose further that
\begin{itemize}
\item[(G4)]$\displaystyle \lim_{s\rightarrow 0^+}\frac{g(s)}{s^{\alpha-1}}=A>0$.
\item[(G5)]$\displaystyle \lim_{s\rightarrow +\infty}\frac{g(s)}{s^{\beta-1}}=B>0$.
\end{itemize}
Then we can obtain the following  asymptotic behavior in the sense of $C_{loc}^{2}(\R^N)$ as well as $H^1(\R^N)$.
\bt\lab{thm:20210814-th1}({\bf The case of $c\rightarrow +\infty$ in the sense of $C_{loc}^{2}(\R^N)$})
Under the assumptions (G1)-(G5), for $c>0$, let $(\lambda_c,u_c)$ be given by Theorem \ref{thm:Main-Thoerem}, then
$$v_c(x):=\lambda_{c}^{\frac{1}{2-\alpha}} u_c\left(\frac{x}{\sqrt{\lambda_c}} \right)\rightarrow U(x)\;\hbox{in}~C_{loc}^{2}(\R^N)~\hbox{as}~c\rightarrow +\infty,$$
where $U$ is the unique radial positive solution of
$$\begin{cases}
-a\Delta U+U=A U^{\alpha-1}~\hbox{in}~\R^N,\\
\lim_{|x|\rightarrow +\infty}U(x)=0.
\end{cases}$$
\et

\bt\lab{thm:20210818-th1}({\bf The case of $c\rightarrow +\infty$ in the sense of $H^1(\R^N)$})
Under the assumptions of (G1)-(G5), let $v_c$ and $U$ be given by Theorem \ref{thm:20210814-th1}, then  we also have that
$$v_c(x)\rightarrow U(x)~\hbox{in}~H^1(\R^N).$$
\et

\bt\lab{thm:20210814-th2}({\bf The case of $c\rightarrow  0^+$ in the sense of $C_{loc}^{2}(\R^N)$})
Under the assumptions (G1)-(G5), for $c>0$, let $(\lambda_c,u_c)$ be given by Theorem \ref{thm:Main-Thoerem}, then
$$v_c(x):=\lambda_{c}^{\frac{1}{2-\beta}} u_c\left(\frac{\|\nabla u_c\|_2}{\sqrt{\lambda_c}}x \right)\rightarrow V(x)\;\hbox{in}~C_{loc}^{2}(\R^N)~\hbox{as}~c\rightarrow 0^+,$$
where $V$ is the unique radial positive solution of
$$\begin{cases}
-b\Delta V+V=B V^{\beta-1}~\hbox{in}~\R^N,\\
\lim_{|x|\rightarrow +\infty}V(x)=0.
\end{cases}$$
\et

\bt\lab{thm:20210818-th2}({\bf The case of $c\rightarrow 0^+$ in the sense of $H^1(\R^N)$})
Under the assumptions of (G1)-(G5), let $v_c$ and $V$ be given by Theorem \ref{thm:20210814-th2}, then  we also have that
$$v_c(x)\rightarrow V(x)~\hbox{in}~H^1(\R^N).$$
\et

\medskip

\section{Preliminaries for the existence}\lab{sec:Preliminaries}
\renewcommand{\theequation}{3.\arabic{equation}}
\bl\lab{lemma:20210807-l1}
If $u\in H^1(\R^N)$ solves
\beq\lab{eq:20210807-e1}
-\left(a+b\int_{\R^N}|\nabla u|^2dx\right)\Delta u+\lambda u=g(u)~\hbox{in}~\R^N,
\eeq
then $u\in \mathcal{P}$, where
$$\mathcal{P}:=\left\{u\in H^1(\R^N): P[u]=0\right\}$$
and
$$P[u]:=a\|\nabla u\|_2^2+b\|\nabla u\|_2^4-N\int_{\R^N}\tilde{G}(u)dx.$$
\el
\bp
Let $u$ be a solution to \eqref{eq:20210807-e1}, we firstly have that
\beq\lab{eq:Nehari}
a\|\nabla u\|_2^2+b\|\nabla u\|_2^4+\lambda \|u\|_2^2=\int_{\R^N}g(u)u dx.
\eeq
Secondly, $u$ satisfies the  so-called Pohozaev identity
\begin{align}\lab{eq:Pohozaev}
&(N-2)\left(a\|\nabla u\|_2^2+b\|\nabla u\|_2^4\right)+N\lambda \|u\|_2^2
-2N\int_{\R^N}G(u)dx=0.
\end{align}
Eliminate the unknown parameter $\lambda$, we obtain that
$$a\|\nabla u\|_2^2+b\|\nabla u\|_2^4-N\int_{\R^N}\tilde{G}(u)dx=0.$$
\ep

\br\lab{remark:20210807-r1}
Recalling the fiber map
$\displaystyle (t\star u)(x):=t^{\frac{N}{2}}u(tx), t\in \R_+,$
which preserves the $L^2$-norm, We define
$$I_{u}(t):=I[t\star u], t>0$$
and remark that
$\displaystyle \frac{1}{t} P[t\star u]=(I_{u})'(t), t>0$. In particular, it holds that $P[u]=(I_{u})'(1)$.
\er

\bl\lab{lemma:20210807-l2}
Let $u\in S_c$. Then: $t\in \R^+$ is a critical point for $I_{u}(t)$ if and only if $t\star u\in \mathcal{P}_{c}$, where
\beq\lab{eq:def-Pc}
\mathcal{P}_{c}:=\mathcal{P}\cap S_c.
\eeq
\el
\bp
It follows by the fact of $\frac{1}{t} P[t\star u]=(I_{u})'(t)$.
\ep

\bl\lab{lemma:20210807-l3}
For any critical point of $I\big|_{\mathcal{P}_{c}}$ , if $(I_{u})''(1)\neq 0$, then there exists some $\lambda\in \R$ such that
$$I'[u]+\lambda u=0~\hbox{in}~\R^N.$$
\el
\bp
Let $u$ be a critical point of $I[u]$ constraint on $\mathcal{P}_{c}$, then there exist $\lambda,\mu\in \R$ such that
\beq\lab{eq:20218507-e2}
I'[u]+\lambda u+\mu P'(u)=0~\hbox{in}~\R^N.
\eeq
We only need to prove that $\mu=0$.
Noting that a function $u$ solves \eqref{eq:20218507-e2} must satisfy the corresponding Pohozave identity
$$(J_u)'(1):=\frac{d}{dt} J[t\star u]\Big|_{t=1}=0,$$
where $J[u]:=I[u]+\frac{1}{2}\lambda \|u\|_2^2+\mu P[u]$ is the corresponding energy functional of \eqref{eq:20218507-e2}.
Observing that
$$J_u(t):=J[t\star u]=I[t\star u]+\frac{1}{2}\lambda \|u\|_2^2+\mu P[t\star u]
=I_{u}(t)+\frac{1}{2}\lambda \|u\|_2^2+\mu t(I_{u})'(t),$$
we have
\begin{align*}
(J_u)'(t)=\frac{d}{dt} J[t\star u]=& (1+\mu)(I_{u})'(t)+\mu t(I_{u})''(t).
\end{align*}
Hence,
\begin{align*}
0=&(J_u)'(1)=(1+\mu)(I_{u})'(1)+\mu (I_{u})''(1)\\
=&(1+\mu)P[u]+\mu ~(I_{u})''(1)=\mu ~(I_{u})''(1).
\end{align*}
So by $(I_{u})''(1)\neq 0$, we have that $\mu=0$.
\ep

\bl\lab{lemma:non-existence}
Let $a,b>0$.
 Assume that there exist some $2<\alpha\leq \beta<2^*$ such that
\beq\lab{eq:AR-condition}
\alpha G(s)\leq g(s)s\leq \beta G(s),~\forall s\in \R.
\eeq
Then equation \eqref{eq:20210807-e1}
has no nontrivial  solution $u\in H^1(\R^N)$ provided $\lambda\leq 0$.
\el
\bp
We argue by contradiction and suppose that there exists some  $0\not\equiv u\in H^1(\R^N)$ solving \eqref{eq:20210807-e1} with $\lambda\leq 0$. By \eqref{eq:Nehari} and \eqref{eq:Pohozaev}, since $\lambda\leq 0$, we have that
$$0\geq \lambda \|u\|_2^2=\int_{\R^N}\left[NG(u)-\frac{N-2}{2}g(u)u\right]dx.$$
For $N\leq 3$, by \eqref{eq:AR-condition}, one can see that $\displaystyle\int_{\R^N}\left[NG(u)-\frac{N-2}{2}g(u)u\right]dx\geq 0$. So we obtain that  $\lambda=0$ and
$$\int_{\R^N}G(u)dx=\int_{\R^N}g(u)udx=\int_{\R^N}\tilde{G}(u)=0.$$
However, even for $\lambda=0$, noting that $u\in \mathcal{P}$, $P[u]=0$, one can also see that
$\|\nabla u\|_2=0$, a contradiction to $0\not\equiv u\in H^1(\R^N)$.
\ep

\bl\lab{lemma:20210807-xl1}
Under the assumptions (G1) and (G2), for any given $c>0$, there exists some $\delta_c>0$ such that
\beq\lab{eq:20210611-e1}
\inf\left\{t>0:\exists u\in S_c~\hbox{with}~\|\nabla u\|_2=1~\hbox{such that}~t\star u\in \mathcal{P}_{c}\right\}\geq \delta_c.
\eeq
That is,
$$\inf\left\{\|\nabla u\|_2, u\in \mathcal{P}_c\right\}\geq \delta_c.$$
\el
\bp
By $t\star u\in \mathcal{P}_c$, we have that $P[t\star u]=0$.  Recalling that $(I_u)'(t)=\frac{1}{t} P[t\star u]$, we get
$$a\|\nabla u\|_2^2+b\|\nabla u\|_2^4 t^2=N\int_{\R^N}\tilde{G}(t^{\frac{N}{2}}u(x))dx t^{-N-2}.$$
By the assumption (G2), for $\|\nabla u\|_2=1$, we have
\beq\lab{eq:20210807-e3}
a<a+bt^2\leq N\left(\frac{1}{2}-\frac{1}{\beta}\right)t^{-N-2}\int_{\R^N} g(t^{\frac{N}{2}}u(x))t^{\frac{N}{2}}u(x) dx.
\eeq
Also by (G2), one can find some suitable $C_1>0$ such that
$$g(s)s\leq C_1\left(|s|^\alpha+|s|^\beta\right), ~\forall s\in \R.$$
Then by Gagliardo-Nierenberg inequality, there exists some $C_2>0$ such that
\beq\lab{eq:20210807-e4}
\|u\|_\alpha^\alpha\leq C_2, \|u\|_\beta^\beta\leq C_2, \forall u\in S_c~\hbox{with}~\|\nabla u\|_2=1.
\eeq
Then by \eqref{eq:20210807-e3} and \eqref{eq:20210807-e4},
$$a<C_1C_2N \left(\frac{1}{2}-\frac{1}{\beta}\right)\left[t^{\frac{N}{2}\alpha -N-2}+t^{\frac{N}{2}\beta -N-2}\right],$$
which implies the existence of lower bound  $\delta_c>0$.
\ep

\bl\lab{lemma:20210807-xl2}
Under the assumptions (G1)-(G3), we have that $(I_u)''(1)<0$ for any $u\in \mathcal{P}_c$ and thus $\mathcal{P}_c$ is a natural constraint of $I\big|_{S_c}$.
\el
\bp
By $(I_u)'(t)=\frac{1}{t} P[t\star u]$, we have that
\begin{align*}
(I_u)''(t)=&a\|\nabla u\|_2^2 +3b\|\nabla u\|_2^4 t^2 +N(N+1)\int_{\R^N}\tilde{G}(t^{\frac{N}{2}}u(x))dx t^{-N-2}\\
&-N\int_{\R^N} \tilde{G}'(t^{\frac{N}{2}}u(x))\frac{N}{2} t^{\frac{N-2}{2}}u(x)dx t^{-N-1}.
\end{align*}
Thus,
\beq\lab{eq:20210807-xe1}
(I_u)''(1)=a\|\nabla u\|_2^2 +3b\|\nabla u\|_2^4+N(N+1)\int_{\R^N}\tilde{G}(u)dx-\frac{N^2}{2}\int_{\R^N}\tilde{G}'(u)udx.
\eeq
Then by the assumption (G3) and $P[u]=0$,
\begin{align*}
(I_u)''(1)\leq &a\|\nabla u\|_2^2 +3b\|\nabla u\|_2^4+N(N+1)\int_{\R^N}\tilde{G}(u)dx
-\frac{N^2}{2}\alpha \int_{\R^N}\tilde{G}(u)dx\\
=& a\|\nabla u\|_2^2 +3b\|\nabla u\|_2^4+[N+1-\frac{N}{2}\alpha] [a\|\nabla u\|_2^2 +b\|\nabla u\|_2^4]\\
=&\left[N+2-\frac{N}{2}\alpha\right] a\|\nabla u\|_2^2+\left[N+4-\frac{N}{2}\beta\right] b\|\nabla u\|_2^4.
\end{align*}
By $2+\frac{8}{N}<\alpha\leq \beta$ and Lemma \ref{lemma:20210807-xl1}, we have that
\beq\lab{eq:20210807-xe2}
(I_u)''(1)\leq \left[N+2-\frac{N}{2}\alpha\right] a\delta_c^2+\left[N+4-\frac{N}{2}\beta\right] b\delta_c^4 <0.
\eeq
Hence, by Lemma \ref{lemma:20210807-l3}, one can see that $\mathcal{P}_c$ is a natural constraint of $I\big|_{S_c}$.
\ep

\br\lab{remark:20210807-xr1}
Let $\{u_n\}\subset \mathcal{P}_{c}$ be such that $I[u_n]$ approaches a possible critical value $M_c$. Noting that $(I_{u_n})''(1)<0$ is an open constraint, then exist sequences $\lambda_n,\mu_n\in \R$ such that
$$I'[u_n]+\lambda_n u_n +\mu_n \mathcal{P}'[u_n]\rightarrow 0.$$
Applying a similar argument as Lemma \ref{lemma:20210807-l3}, we obtain that
\beq\lab{eq:20210807-xe3}
\mu_n (I_{u_n})''(1)\rightarrow 0.
\eeq
By \eqref{eq:20210807-xe2}, one can see that
$$\mu_n \left\{\left[N+2-\frac{N}{2}\alpha\right] a\delta_c^2+\left[N+4-\frac{N}{2}\beta\right] b\delta_c^4\right\}\rightarrow 0.$$
However, by Lemma \ref{lemma:20210807-xl1},  we have that $\|\nabla u\|_2^2\geq \delta_c^2>0$ for all $u\in \mathcal{P}_c$. Hence, we obtain that $\mu_n\rightarrow 0$.
And thus if furthermore $\{u_n\}$ is bounded in $H^1(\R^N)$, then we have that
$$I'[u_n]+\lambda_n u_n\rightarrow 0\;\hbox{in}\;H^{-1}(\R^N).$$
\er

\bc\lab{cro:20210807-c1}
Under the assumptions (G1)-(G3), for any $u\in H^1(\R^N)\backslash\{0\}$, there exists an unique $t_u>0$ such that $t_u\star u\in \mathcal{P}$. Furthermore,
\beq\lab{eq:20210807-xe4}
I[t_u\star u]=\max_{t>0}I[t\star u].
\eeq
\ec
\bp
Let $c:=\|u\|_2^2>0$.
We remark that under the assumptions (G1) and (G2), for any $t\geq 0$ and $s\in \R$, it holds that
\beq\lab{eq:20210807-xe4}
\begin{cases}
t^\beta G(s)\leq G(ts)\leq t^\alpha G(s),\quad & \hbox{if}~t\leq 1,\\
t^\alpha G(s)\leq G(ts)\leq t^\beta G(s),\quad & \hbox{if}~t\geq 1.
\end{cases}
\eeq
Then we have that
\beq\lab{eq:20210807-xe5}
\frac{\alpha-2}{\beta-2} \min\{t^\alpha, t^\beta\}  \tilde{G}(s)\leq \tilde{G}(ts)\leq
\frac{\beta-2}{\alpha-2} \max\{t^\alpha, t^\beta\}  \tilde{G}(s),
\eeq
and thus by $\alpha>2+\frac{8}{N}$,
\beq\lab{eq:20210807-xe6}
\int_{\R^N}\tilde{G}(t^{\frac{N}{2}}u)dx t^{-N}=o(t^4)~\hbox{as}~t\rightarrow 0^+.
\eeq
Hence, by
$$
P[t\star u]=a\|\nabla u\|_2^2 t^2+b\|\nabla u\|_2^4 t^4 -N\int_{\R^N}\tilde{G}(t^{\frac{N}{2}}u)dx t^{-N},
$$
one can see that $P[t\star u]>0$ for $t>0$ small enough.
Recalling that $(I_u)'(t)=\frac{1}{t}P[t\star u]$, we have that
$(I_u)'(t)>0$ for $t>0$ small enough. Hence, there exists some $t_0>0$ such that $I_u(t)$ increases in $t\in (0,t_0)$.

On the other hand, by $\displaystyle \int_{\R^N}G(u)dx>0$ and $\alpha>2+\frac{8}{N}$, we can also have that
\beq\lab{eq:20210807-xe7}
\int_{\R^N}G(t^{\frac{N}{2}}u)dx t^{-N-4}\geq \left(\int_{\R^N}G(u)dx\right) t^{\frac{N}{2}\alpha -N-4} \rightarrow +\infty~\hbox{as}~t\rightarrow +\infty.
\eeq
 Thus
\beq\lab{eq:20210807-xe8}
\lim_{t\rightarrow +\infty} I_{u}(t)=\lim_{t\rightarrow +\infty} t^4\left\{\frac{a}{2}\|\nabla u\|_2^2 t^{-2}+\frac{b}{4}\|\nabla u\|_2^4-\int_{\R^N}G(t^{\frac{N}{2}}u)dx t^{-N-4}\right\}=-\infty.
\eeq
So there must exists some $t_1>t_0$ such that
$$I_u(t_1)=\max_{t>0}I_u(t).$$
Hence, $(I_u)'(t_1)=0$ and it follows Lemma \ref{lemma:20210807-l2} that $t_1\star u\in \mathcal{P}$.
Suppose that there exists another $t_2>0$ such that $t_2\star u\in \mathcal{P}$. Then by Lemma \ref{lemma:20210807-xl2}, we have that both $t_1$ and $t_2$ are strict local maximum of $I_u(t)$. Without loss of generality, we assume that $t_1<t_2$. Then there exists some $t_3\in (t_1, t_2)$ such that
$$I_u(t_3)=\min_{t\in [t_1,t_2]}I_u(t).$$
That is, $t_3$ is a local minimum of $I_u(t)$.
So we have that $(I_u)'(t_3)=0$ and thus $t_3\star u\in \mathcal{P}$ with $(I_{t_3\star u})''(1)=(I_{u})''(t_3)\geq 0$, a contradiction to Lemma  \ref{lemma:20210807-xl2}.
\ep

\bc\lab{cro:20210807-wc1}
Under the assumptions (G1)-(G3), for any $u\in H^1(\R^N)\backslash\{0\}$, let $t_u$ be given by Corollary \ref{cro:20210807-c1}, then we have that
$$t_u=(>,<)1\Leftrightarrow (I_u)'(1)=(>,<)0 \Leftrightarrow P[u]=(>,<)0.$$
\ec
\bp
By Corollary \ref{cro:20210807-c1}, we have that
$$I_u(t_u)=\max_{t>0}I_u(t).$$
 Furthermore,
 $$(I_u)'(t)>0\;\hbox{for}~0<t<t_u~\hbox{and}~(I_u)'(t)<0\;\hbox{for}~t>t_u.$$
 On the other hand, we recall that $P[t\star u]=t(I_u)'(t)$.
 Hence, the conclusion holds.
\ep

\br\lab{remark:Mc-positive}
Under the hypotheses (G1)-(G3), for any $u\in \mathcal{P}_{c}$,
one can see that $I_u(t)\rightarrow 0$ as $t\rightarrow 0^+$ and $I_u(t)\rightarrow -\infty$ as $t\rightarrow +\infty$. By Corollary \ref{cro:20210807-c1}, we have that
$$I[u]=\max_{t>0}I[t\star u]>0.$$
\er

\bl\lab{lemma:20210807-xl3}
Under the assumptions (G1)-(G3), $I\big|_{\mathcal{P}_{c}}$ is coercive , i.e.,
$$\lim_{u\in \mathcal{P}_{c}, \|\nabla u\|_2\rightarrow \infty}I[u]=+\infty.$$
\el
\bp
For $u\in \mathcal{P}_c$, by (G2), we have that
\begin{align*}
a\|\nabla u\|_2^2+b\|\nabla u\|_2^4=&N\int_{\R^N}\tilde{G}(u)dx
\geq N\frac{\alpha-2}{2}\int_{\R^N}G(u)dx.
\end{align*}
Thus,
\begin{align}\lab{eq:20210807-xe9}
I[u]=&\frac{a}{2}\|\nabla u\|_2^2+\frac{b}{4}\|\nabla u\|_2^4-\int_{\R^N}G(u)dx\nonumber\\
\geq&\frac{a}{2}\|\nabla u\|_2^2+\frac{b}{4}\|\nabla u\|_2^4 - \frac{2}{N(\alpha-2)}\left(a\|\nabla u\|_2^2+b\|\nabla u\|_2^4\right)\nonumber\\
=&\frac{N(\alpha-2)-4}{2N(\alpha-2)}a\|\nabla u\|_2^2+\frac{N(\alpha-2)-8}{4N(\alpha-2)} b\|\nabla u\|_2^4.
\end{align}
The conclusion follows that  $\frac{N(\alpha-2)-4}{2N(\alpha-2)}>0,\frac{N(\alpha-2)-8}{4N(\alpha-2)}>0$ due to $\alpha>2+\frac{8}{N}$.
\ep

For given $c>0$, define
\beq\lab{eq:20210807-xe10}
M_c:=\inf_{u\in \mathcal{P}_c}I[u]=\inf_{u\in S_c}\max_{t>0}I[t\star u].
\eeq
Since a solution $u$ to \eqref{eq:Main-Problem} with \eqref{eq:20210816-constraint} must belong to $\mathcal{P}_c$, one can see that $u$ is indeed a least energy solution (ground state solution) if $u$ attains $M_c$.

\bc\lab{cro:20210807-xc2}
Under the assumptions (G1)-(G3), for any $c>0$, we have $M_c>0$.
\ec
\bp
It follows by the formula \eqref{eq:20210807-xe9} and Lemma \ref{lemma:20210807-xl1}.
\ep

For any $u\in H^1(\R^N)$, let $u^*$ be the symmetric decreasing rearrangement of $u$. Then
we have that $u^*=|u|^*$. Under the assumption (G1), without loss of generality, we may suppose that $u$ is nonnegative. Then one can see that
\begin{align*}
\int_{\R^N}G(u)dx=&\int_{\R^N}\left(\int_{0}^{u(x)}g(s)ds\right) dx\\
=&\int_0^\infty g(s)\left|\left\{x: u(x)>s\right\}\right|ds\\
=&\int_0^\infty g(s)\left|\left\{x: u^*(x)>s\right\}\right|ds\\
=&\int_{\R^N}G(u^*)dx.
\end{align*}
Recalling the P\'olya-Szeg\"o inequality, we have
\beq\lab{eq:Polya-Szego}
\int_{\R^N}|\nabla u^*|^2 dx\leq \int_{\R^N}|\nabla u|^2 dx.
\eeq
Hence,
\beq\lab{eq:20210807-we1}
I[u^*]\leq I[u].
\eeq

Let $H^{1}_{rad}(\R^N)$ be the radial subspace of $H^1(\R^N)$. Set
\beq
S_{c}^{rad}:=S_c \cap H^{1}_{rad}(\R^N), \mathcal{P}^{rad}:=\mathcal{P}\cap H^{1}_{rad}(\R^N), \mathcal{P}_{c}^{rad}:=\mathcal{P}_c\cap H^{1}_{rad}(\R^N).
\eeq

Define
\beq\lab{eq:Def-Mc-rad}
M_{c}^{rad}:=\inf_{u\in S_{c}^{rad}}\max_{t>0}I[t\star u].
\eeq
Then one can see that it also holds that
\beq\lab{eq:20210807-wbe1}
M_{c}^{rad}=\inf_{u\in \mathcal{P}_{c}^{rad}}I[u].
\eeq
Furthermore, we have:
\bl\lab{lemma:20210807-wbl1}
$$M_{c}^{rad}=M_c.$$
\el
\bp
Since $S_{c}^{rad}\subset S_c$, it is trivial that $M_{c}^{rad}\geq M_c$. On the other hand, for any $s>0$,
\begin{align*}
&|\{x:t\star u^*(x)>s\}|=|\{x:t^{\frac{N}{2}}u^*(tx)>s\}|=|\{y: t^{\frac{N}{2}}u^*(y)>s\}|t^{-N}\\
=&t^{-N}|\{u^*(y)>t^{-\frac{N}{2}}s\}|=t^{-N}|\{u(y)>t^{-\frac{N}{2}}s\}|
=t^{-N}|\{t^{\frac{N}{2}}u(y)>s\}|\\
=&|\{t^{\frac{N}{2}}u(tx)>s\}|
=|\{(t\star u)(x)>s\}|
=|\{(t\star u)^*(x)>s\}|,
\end{align*}
one can see that $$t\star u^*=(t\star u)^*, \forall t\in \R^+.$$
Then for any $u\in \mathcal{P}_c$, we have
\begin{align*}
I[t\star u^*]=&I[(t\star u)^*]\leq I[t\star u]\\
\leq&\max_{s>0}I[s\star u]=I[u], \forall t\in \R^+.
\end{align*}
Then by the arbitrary of $u\in \mathcal{P}_c$, we obtain that $M_{c}^{rad}\leq M_c$. And thus $M_{c}^{rad}=M_c$.
\ep

\section{Proof of Theorem \ref{thm:Main-Thoerem}}\lab{Proof-th1}
\renewcommand{\theequation}{4.\arabic{equation}}
\bp
Let $\{u_n\}\subset \mathcal{P}_{c}^{rad}$ be such that $I[u_n]\rightarrow M_c>0$. By Lemma \ref{lemma:20210807-xl3}, we see that $\{u_n\}$ is bounded in $H^1(\R^N)$. Up to a subsequence, we may assume that $u_n\rightharpoonup u$ in $H^1(\R^N)$. Then for $N=2,3$, by the compact embedding $H_{rad}^{1}(\R^N)\hookrightarrow\hookrightarrow L^p(\R^N)$ for $2<p<2^*$, we have that
\beq\lab{eq:20210807-wbe2}
\int_{\R^N}G(u_n)dx\rightarrow \int_{\R^N}G(u)dx.
\eeq
For the case of $N=1$, we may suppose further that $u_n=u_n^*, \forall n\in \N$. Then
\eqref{eq:20210807-wbe2} also holds (see \cite[Proposition 1.7.1]{Cazenave2003}).

We claim that $u\neq 0$. If not, $\int_{\R^N}\tilde{G}(u_n)dx=o(1)$ and thus $\{u_n\}\subset \mathcal{P}_{c}^{rad}$ implies that
$$a\|\nabla u_n\|_2^2+b\|\nabla u_n\|_2^4=o(1),$$
a contradiction to Lemma \ref{lemma:20210807-xl1}.

Since $\{u_n\}$ is bounded in $H^1(\R^N)$, it is easy to see that
$$\lambda_n:=-\frac{1}{c} \langle I'[u_n], u_n\rangle$$
is a bounded sequence. In particular,
\begin{align}\lab{eq:20210810-be1}
\lambda_n c=&\lambda_n\|u_n\|_2^2=-\langle I'[u_n], u_n\rangle\nonumber\\
=&\int_{\R^N}g(u_n)u_ndx-a\|\nabla u_n\|_2^2-b\|\nabla u_n\|_2^4\nonumber\\
=&\int_{\R^N}g(u_n)u_ndx-N\int_{\R^N}\tilde{G}(u_n)dx\nonumber\\
=&N\int_{\R^N}G(u_n)dx-\frac{N-2}{2}\int_{\R^N}g(u_n)u_ndx\nonumber\\
\geq&\left[N-\frac{(N-2)\beta}{2}\right]\int_{\R^N}G(u_n)dx.
\end{align}
On the other hand,
\begin{align*}
a\|\nabla u_n\|_2^2+b\|\nabla u_n\|_2^4=N\int_{\R^N}\tilde{G}(u_n)dx\leq \frac{(\beta-2)N}{2}\int_{\R^N}G(u_n)dx.
\end{align*}
Since $2+\frac{8}{N}<\alpha\leq \beta<2^*$, by Lemma \ref{lemma:20210807-xl1}, one can find some $\eta_c>0$ such that
$$\int_{\R^N}G(u_n)dx\geq \eta_c, \forall n\in \N.$$
And thus there exists some suitable $\sigma_c>0$ such that
$$\lambda_n\geq \sigma_c, \forall n\in \N.$$

Without loss of generality, we may assume that $\lambda_n\rightarrow \lambda>0$.
Suppose that up to a subsequence, $\|\nabla u_n\|_2^2\rightarrow \Lambda \geq 0$. Then
one can see that $u\in H^{1}_{rad}(\R^N)$ solves
\beq\lab{eq:20210818-xbe1}
-\left(a+b\Lambda \right) \Delta u+\lambda u=g(u)~\hbox{in}~\R^N.
\eeq
Then we also have that $u\in \mathcal{P}^{rad}$, and thus
\begin{align*}
a\|\nabla u\|_2^2+b\Lambda \|\nabla u\|_2^2=&\int_{\R^N}\tilde{G}(u)dx\\
=&\lim_{n\rightarrow \infty} \int_{\R^N}\tilde{G}(u_n)dx\\
=&\lim_{n\rightarrow \infty} \left(a\|\nabla u_n\|_2^2+b\|\nabla u_n\|_2^4\right)\\
=&a\Lambda +b\Lambda^2.
\end{align*}
So $(a+b\Lambda)(\|\nabla u\|_2^2-\Lambda)=0$. By $a>0,b>0,\Lambda\geq 0$, we see that $\|\nabla u\|_2^2=\Lambda$, which implies that $u_n\rightarrow u$ in $D_{0}^{1,2}(\R^N)$.
So by \eqref{eq:20210818-xbe1}, we see that $u$ solves
\beq\lab{eq:20210818-xbe2}
-\left(a+b\int_{\R^N}|\nabla u|^2dx \right) \Delta u+\lambda u=g(u)~\hbox{in}~\R^N.
\eeq
On the other hand, we also have that
\begin{align*}
a\|\nabla u\|_2^2+b\|\nabla u\|_2^4+\lambda \|u\|_2^2=&\int_{\R^N}g(u)udx\\
=&\int_{\R^N}g(u_n)u_ndx+o(1)\\
=&a\|\nabla u_n\|_2^2+b\|\nabla u_n\|_2^4+\lambda_n \|u_n\|_2^2+o(1),
\end{align*}
which implies that
$$\lambda(c-\|u\|_2^2)=0.$$
We obtain that $u\in S_c$ and thus $u\in \mathcal{P}_c$.
Hence,
$$M_c\leq I[u]=\lim_{n\rightarrow \infty}I[u_n]=M_c.$$
That is, $(\lambda, u)$ is a  ground state normalized solution of Kirchhoff equation \eqref{eq:Main-Problem}.
\ep

\section{Properties of $M_c$ and $\lambda_c$}\lab{sec:Properties-Mc}
\renewcommand{\theequation}{5.\arabic{equation}}
In this section, we shall study some properties of $M_c$ and $\lambda_c$, such as the continuity and the limit behavior as $c\rightarrow 0^+$ as well as $c\rightarrow +\infty$, which is very important for us to study the asymptotic behavior in Section \ref{sec:aysmptotic-behavior}.

Now, for $c>0$, let $(\lambda_c, u_c)$ be the solution given by Theorem \ref{thm:Main-Thoerem}. That is, $\lambda_c>0$ and $u_c\in S_{c}^{rad}$ satisfy
\beq\lab{eq:20210810-e1}
-\left(a+b\int_{\R^N}|\nabla u_c|^2dx\right)\Delta u_c+\lambda_c u_c=g(u_c)~\hbox{in}~\R^N,
\eeq
and
\beq\lab{eq:20210810-e2}
I[u_c]=\frac{a}{2}\|\nabla u_c\|_2^2+\frac{b}{4}\|\nabla u_c\|_2^4-\int_{\R^N}G(u_c)dx=M_c.
\eeq

Then we have the following results:
\bl\lab{lemma:20210810-l1}
$M_c$ is continuous with respect to $c\in (0,+\infty)$.
\el
\bp
Noting that $M_c$ is attained by some $u_c\in H_{rad}^{1}(\R^N)$, which is a symmetric decreasing function. Fix $c>0$, for any $\{c_n\}\subset \R_+$ with $c_n\rightarrow c$ as $n\rightarrow +\infty$, we simply write $(\lambda_{c_n},u_{c_n})$ by $(\lambda_n,u_n)$.
Without loss of generality, we may assume that $c_n\in (\frac{c}{2}, 2c), \forall n\in \N$.
Noting that
\begin{align}\lab{eq:20210810-e7}
&\left[\frac{1}{2}-\frac{2}{(\alpha-2)N}\right] a\|\nabla u_c\|_2^2+\left[\frac{1}{4}-\frac{2}{(\alpha-2)N}\right] b\|\nabla u_c\|_2^4\nonumber\\
\leq &M_c=I[u_c]\nonumber\\
\leq &\left[\frac{1}{2}-\frac{2}{(\beta-2)N}\right] a\|\nabla u_c\|_2^2+\left[\frac{1}{4}-\frac{2}{(\beta-2)N}\right] b\|\nabla u_c\|_2^4.
\end{align}
Let $u_{\frac{c}{2}}$ attain $M_{\frac{c}{2}}$. For any $\kappa\in (1,4)$, $\kappa u_{\frac{c}{2}}\in S_{\kappa^2 \frac{c}{2}}$. Recalling that there exists an unique $t_\kappa>0$ such that $t_\kappa\star (\kappa u_{\frac{c}{2}})\in \mathcal{P}_{\kappa^2 \frac{c}{2}}$, and thus
\beq
a \|\kappa\nabla u_{\frac{c}{2}}\|_2^2 t_{\kappa}^{2} +b \|\kappa\nabla u_{\frac{c}{2}}\|_2^4 t_{\kappa}^{4} =N\int_{\R^N}\tilde{G}(t_\kappa\star (\kappa u_{\frac{c}{2}})).
\eeq
In particular, by \eqref{eq:20210807-xe5}, we have
\begin{align*}
a \|\kappa\nabla  u_{\frac{c}{2}}\|_2^2 t_{\kappa}^{-2} +b \|\kappa\nabla  u_{\frac{c}{2}}\|_2^4 =&Nt_{\kappa}^{-4}\int_{\R^N}\tilde{G}(t_\kappa\star (\kappa u_{\frac{c}{2}}))\\
\geq & C t_{\kappa}^{-4} \int_{\R^N}\tilde{G}(t_{\kappa}^{\frac{N}{2}} \kappa u_{\frac{c}{2}}(t_{\kappa} x))dx\\
\geq &C t_{\kappa}^{-4} \min\{(t_{\kappa}^{\frac{N}{2}} \kappa)^\alpha, (t_{\kappa}^{\frac{N}{2}} \kappa )^\beta \} \int_{\R^N}\tilde{G}(u_{\frac{c}{2}}(t_{\kappa} x))dx\\
=&C t_{\kappa}^{-4-N} \min\{(t_{\kappa}^{\frac{N}{2}} \kappa)^\alpha, (t_{\kappa}^{\frac{N}{2}} \kappa )^\beta \} \int_{\R^N}\tilde{G} (u_{\frac{c}{2}}(y))dy\\
\geq &C \min\{t_{\kappa}^{\frac{N}{2}\alpha -N-4}, t_{\kappa}^{\frac{N}{2}\beta -N-4}\},
\end{align*}
which implies that $t_\kappa$ is bounded since $\frac{N}{2}\beta -N-4\geq \frac{N}{2}\alpha -N-4>0$.
Hence, by $\displaystyle M_{\kappa^2 \frac{c}{2}}\leq I[t_\kappa\star (\kappa u_{\frac{c}{2}})]$, we see  that $\{M_{\kappa^2 \frac{c}{2}}: \kappa\in (1,4)\}$ is bounded. So by the left hand side of \eqref{eq:20210810-e7}, we obtain that $\{u_n\}$ is bounded in $H^1(\R^N)$.

Similar to the proof of Theorem \ref{thm:Main-Thoerem}, we may assume that $\lambda_n\rightarrow \lambda>0$ and $u_n\rightharpoonup u$ in $H^1(\R^N)$. Here $u\in H^1(\R^N)$ satisfies $u\in S_c$ and solves
$$-\left(a+b\int_{\R^N}|\nabla u|^2 dx\right)\Delta u+\lambda u=g(u)~\hbox{in}~\R^N.$$
Hence,
\beq\lab{eq:20210810-e8b}
\lim_{n\rightarrow \infty}M_{c_n}=\lim_{n\rightarrow \infty}I[u_n]=I[u]\geq M_c.
\eeq
If $I[u]\neq M_c$, then there exists some $\delta>0$  such that
\beq\lab{eq:20210810-e9}
M_{c_n}\geq M_c+\delta~\hbox{for $n$ large enough}.
\eeq
We remark that $\sqrt{\theta}u_c\in S_{\theta c}$. Let $t_{\theta}>0$ be the unique number such that $t_\theta \star (\sqrt{\theta}u_c)= \sqrt{\theta}(t_\theta \star u_c)\in \mathcal{P}_{\theta c}$. That is,
\beq\lab{eq:20210810-e10}
0=a\|\nabla u_c\|_2^2 \theta t_\theta^2 +b\|\nabla u_c\|_2^4 \theta^2 t_\theta^4 -\int_{\R^N}\tilde{G}(\sqrt{\theta}(t_\theta \star u_c)).
\eeq
Recalling Lemma \ref{lemma:20210807-xl1}, one can prove that $t_\theta$ is bounded away from $0$ for $\theta$ close to $1$. And then the uniqueness implies that $t_\theta\rightarrow 1$ as $\theta\rightarrow 1$.
So that, for $\theta$ close to $1$ enough, we have that
\beq\lab{eq:20210810-e11}
M_{\theta c}\leq I[\sqrt{\theta}\star (t_\theta\star u_c)]\rightarrow I[u_c]=M_c.
\eeq
Then, there exists some $N_0\in \N$, such that
\beq\lab{eq:20210810-e12}
M_{c_n}\leq M_c+\frac{\delta}{2}\;\hbox{for $n\geq N_0$},
\eeq
a contradiction to \eqref{eq:20210810-e8b}.
Hence, we prove that $I[u]=M_c$ and thus
$$\lim_{n\rightarrow \infty}M_{c_n}=M_c.$$
\ep

\bl\lab{lemma:20210810-l2}
$$\lim_{c\rightarrow 0^+}M_c=+\infty~\hbox{and}~\lim_{c\rightarrow +\infty}M_c=0.$$
\el
\bp
By \eqref{eq:20210807-xe4} and the assumption (G2), combing with the Gagliardo-Nierenberg inequality, there exists some $C=C(\alpha,\beta,N)$ such that
\begin{align*}
a\|\nabla u_c\|_2^2+b\|\nabla u_c\|_2^4=&N\int_{\R^N}\tilde{G}(u_c)dx\\
\leq&\frac{(\beta-2)N}{2}\int_{\R^N}G(u_c)dx\\
\leq& C \left(\|u_c\|_\alpha^\alpha+\|u_c\|_\beta^\beta\right)\\
\leq&C \left(\|\nabla u_c\|_{2}^{\frac{N(\alpha-2)}{2}} c^{\frac{2\alpha-N(\alpha-2)}{4}}+\|\nabla u_c\|_{2}^{\frac{N(\beta-2)}{2}} c^{\frac{2\beta-N(\beta-2)}{4}}\right).
\end{align*}
Thus,
\beq\lab{eq:20210810-e3}
a<a+b\|\nabla u_c\|_2^2\leq C\left(\|\nabla u_c\|_{2}^{\frac{N(\alpha-2)-4}{2}} c^{\frac{2\alpha-N(\alpha-2)}{4}}+\|\nabla u_c\|_{2}^{\frac{N(\beta-2)-4}{2}} c^{\frac{2\beta-N(\beta-2)}{4}}\right).
\eeq
Noting that $\frac{2\alpha-N(\alpha-2)}{4}$ and  $\frac{2\beta-N(\beta-2)}{4}$ are positive, one can see that
\beq\lab{eq:20210810-e4}
\|\nabla u_c\|_{2}^{\frac{N(\alpha-2)-4}{2}}, \|\nabla u_c\|_{2}^{\frac{N(\beta-2)-4}{2}}\rightarrow +\infty~\hbox{as}~c\rightarrow 0^+.
\eeq
So by $\frac{N(\alpha-2)-4}{2},\frac{N(\beta-2)-4}{2}>0$, we obtain that
\beq\lab{eq:20210810-e5}
\|\nabla u_c\|_2\rightarrow +\infty~\hbox{as}~c\rightarrow 0^+.
\eeq
Recalling \eqref{eq:20210807-xe9}, we obtain that
\beq\lab{eq:20210810-e6}
M_c=I[u_c]\geq \frac{N(\alpha-2)-4}{2N(\alpha-2)}a\|\nabla u_c\|_2^2+\frac{N(\alpha-2)-8}{4N(\alpha-2)} b\|\nabla u_c\|_2^4\rightarrow +\infty~\hbox{as}~c\rightarrow 0^+.
\eeq

On the other hand, we fix some positive $\omega\in H^1(\R^N)$ with $\|\omega\|_2^2=1$. For $c>0$, one has that $u=\sqrt{c}\omega\in S_c$. Let $t_c>0$ be the unique number such that $t_c\star u\in \mathcal{P}_c$.
We claim that $ct_c^2\rightarrow 0$ as $c\rightarrow +\infty$. If not, there exists a sequence $\{c_n\}$ with $c_n\rightarrow +\infty $ as $n\rightarrow +\infty$ and $c_n t_{c_n}^{2}\geq \sigma>0, \forall n\in \N$.
Recalling \eqref{eq:20210807-xe5}, we have that
\begin{align*}
&c_{n}^{-2}t_{c_n}^{-4}\int_{\R^N}\tilde{G}(\sqrt{c_n}(t_{c_n}\star \omega))dx\\
=&c_{n}^{-2}t_{c_n}^{-4}\int_{\R^N}\tilde{G}(\sqrt{c_n}t_{c_n}^{\frac{N}{2}}\omega(y)) t_{c_n}^{-N}dy\\
\geq &c_{n}^{-2} t_{c_n}^{-4-N} \frac{\alpha-2}{\beta-2}\min\left\{\left(\sqrt{c_n}t_{c_n}^{\frac{N}{2}}\right)^\alpha, \left(\sqrt{c_n}t_{c_n}^{\frac{N}{2}}\right)^\beta\right\} \int_{\R^N}\tilde{G}(\omega)dx.
\end{align*}
Noting that
\begin{align*}
&c_{n}^{-2} t_{c_n}^{-4-N} \left(\sqrt{c_n}t_{c_n}^{\frac{N}{2}}\right)^\alpha
=\left(c_nt_{c_n}^{2}\right)^{\frac{1}{2}(\frac{N\alpha}{2}-N-4)} c_{n}^{\frac{2N-(N-2)\alpha}{4}}\\
\geq &\sigma^{\frac{1}{2}(\frac{N\alpha}{2}-N-4)} c_{n}^{\frac{2N-(N-2)\alpha}{4}}
\rightarrow  +\infty~\hbox{as}~n\rightarrow +\infty
\end{align*}
and similarly that
$$c_{n}^{-2} t_{c_n}^{-4-N} \left(\sqrt{c_n}t_{c_n}^{\frac{N}{2}}\right)^\beta \rightarrow  +\infty~\hbox{as}~n\rightarrow +\infty,$$
we have
$$c_{n}^{-2}t_{c_n}^{-4}\int_{\R^N}\tilde{G}(\sqrt{c_n}(t_{c_n}\star \omega))dx\rightarrow +\infty ~\hbox{as}~n\rightarrow +\infty.$$
Then we have that
\begin{align*}
0=&P[t_{c_n}\star u]=a\|\nabla u\|_2^2 t_{c_n}^{2}+b\|\nabla u\|_2^4t_{c_n}^{4}-\int_{\R^N}\tilde{G}(t_{c_n}\star u)dx\\
=&a\|\nabla \omega\|_2^2 c_n t_{c_n}^{2} +b\|\nabla \omega\|_2^4 c_n^2t_{c_n}^{4}-\int_{\R^N}\tilde{G}(\sqrt{c}(t_{c_n}\star \omega))dx\\
=& c_n^2t_{c_n}^{4}\left[ a\|\nabla \omega\|_2^2 c_{n}^{-1} t_{c_n}^{-2}+b\|\nabla \omega\|_2^4 -c_{n}^{-2}t_{c_n}^{-4}\int_{\R^N}\tilde{G}(\sqrt{c_n}(t_{c_n}\star \omega))dx\right]\\
<&0\;\hbox{for $n$ large},
\end{align*}
a contradiction.  Hence, the claim $ct_c^2\rightarrow 0$ as $c\rightarrow +\infty$ is proved.

Hence, by $v:=t_c\star u=\sqrt{c} (t_c\star \omega)\in \mathcal{P}_c$, we obtain that
\begin{align*}
M_c\leq I[v]=&\frac{a}{2}\|\nabla v\|_2^2+\frac{b}{4}\|\nabla v\|_2^4-\int_{\R^N}G(v)dx\\
\leq&\frac{a}{2}\|\nabla v\|_2^2+\frac{b}{4}\|\nabla v\|_2^4 -\frac{2}{\beta-2}\int_{\R^N}\tilde{G}(v)dx\\
=&\frac{a}{2}\|\nabla v\|_2^2+\frac{b}{4}\|\nabla v\|_2^4-\frac{2}{(\beta-2)N}\left[a\|\nabla v\|_2^2+b\|\nabla v\|_2^4\right]\\
=&\left[\frac{1}{2}-\frac{2}{(\beta-2)N}\right]a\|\nabla v\|_2^2+\left[\frac{1}{4}-\frac{2}{(\beta-2)N}\right]b\|\nabla v\|_2^4\\
=&\left[\frac{1}{2}-\frac{2}{(\beta-2)N}\right]a\|\nabla \omega\|_2^2 ct_c^2 +\left[\frac{1}{4}-\frac{2}{(\beta-2)N}\right]b\|\nabla \omega\|_2^4 (ct_c^2)^2\\
\rightarrow &0~\hbox{as}~c\rightarrow +\infty.
\end{align*}
\ep

\bc\lab{cor:20210810-c1}
$$\lim_{c\rightarrow 0^+}\lambda_c c=+\infty~\hbox{and}~\lim_{c\rightarrow +\infty}\lambda_c c=0.$$
\ec

\bp
By Lemma \ref{lemma:20210810-l2} and the formula \eqref{eq:20210810-e7},  we see that
\beq\lab{eq:20210810-e8}
\lim_{c\rightarrow +\infty}\|\nabla u_c\|_2=0~\hbox{and}~\lim_{c\rightarrow 0^+}\|\nabla u_c\|_2=+\infty.
\eeq
Furthermore, by $P[u_c]=0$, we have that
$$\int_{\R^N}\tilde{G}(u_c)dx=\frac{1}{N}\left[a\|\nabla u_c\|_2^2+b\|\nabla u_c\|_2^4\right].$$
Hence, we also have
\beq\lab{eq:20210810-e9}
\lim_{c\rightarrow +\infty}\int_{\R^N}\tilde{G}(u_c)dx=0~\hbox{and}~\lim_{c\rightarrow 0^+}\int_{\R^N}\tilde{G}(u_c)dx=+\infty.
\eeq
Similar to the formula \eqref{eq:20210810-be1}, we have that
$$\lambda_c c=\lambda_c \|u_c\|_2^2=\int_{\R^N}\left[NG(u_c)-\frac{N-2}{2}g(u_c)u_c\right]dx.$$
Under the assumption (G2), we can find some $C_1\geq C_2>0$ such that
\beq\lab{eq:20210819-we1}
C_2\int_{\R^N}\tilde{G}(u_c)dx\leq \lambda_c c\leq C_1\int_{\R^N}\tilde{G}(u_c)dx.
\eeq
Hence, by \eqref{eq:20210810-e9}, we obtain that
$$\lambda_c c\rightarrow +\infty~\hbox{as}~c\rightarrow 0^+~\hbox{and}~\lambda_c c\rightarrow 0~\hbox{as}~c\rightarrow +\infty.$$
\ep

\br\lab{remark:20210812-r1}
We say the two quantities $p(c)$ and $q(c)$ are comparable, if there exist some $C_1\geq C_2>0$ independent of $c$  such that
$$C_2 q(c)\leq p(c)\leq C_1q(c).$$
Then by the proofs in Section \ref{sec:Properties-Mc}, it is not hard to see that any two elements of $$\left\{M_c, \lambda_cc, \|\nabla u_c\|_2^2+\|\nabla u_c\|_2^4, \int_{\R^N}G(u_c)dx, \int_{\R^N}\tilde{G}(u_c)dx, \int_{\R^N}g(u_c)u_cdx\right\}$$
are comparable.
\er

\section{Asymptotic behavior of $u_c$}\lab{sec:aysmptotic-behavior}
\renewcommand{\theequation}{6.\arabic{equation}}
\subsection{The case of $c\rightarrow +\infty$}
Let $c_n\rightarrow +\infty$, by Corollary \ref{cor:20210810-c1}, we have that  $\lambda_{c_n}\rightarrow 0$. Write $(\lambda_{c_n}, u_{c_n})$ simply by $(\lambda_n, u_n)$.
Put
\beq\lab{eq:20210810-we1}
d_n:=\|\nabla u_n\|_2^2.
\eeq
Recalling \eqref{eq:20210810-e8}, we have that $d_n\rightarrow 0$ as $n\rightarrow +\infty$.

\bl\lab{lemma:20210813-l1}
Under the assumptions of (G1)-(G5), let $(\lambda_c,u_c)$ be the solution given by Theorem \ref{thm:Main-Thoerem}. Then
$$\limsup_{c\rightarrow +\infty} u_c(0)<+\infty.$$
\el
\bp
Suppose that there exists a sequence $\{c_n\}\rightarrow +\infty$ such that
$$M_n:=u_n(0)=\max_{x\in \R^N}u_n(x)\rightarrow +\infty.$$
Now we perform a rescaling, setting $x=\frac{y}{M_{n}^{\frac{\beta-2}{2}}}$ and defining
$$U_n(y)=\frac{u_n\left(\frac{y}{M_{n}^{\frac{\beta-2}{2}}}\right)}{M_n}, y\in \R^N.$$
Then
$$1=\max_{y\in \R^N}U_n(y)$$
and
\beq\lab{eq:20210813-xe1}
-\Delta U_n(y)=\frac{1}{a+bd_n}\left[\frac{g(M_n U_n(y))}{M_{n}^{\beta-1}}-\frac{\lambda_n}{M_{n}^{\beta-2}}U_n(y)\right].
\eeq
Noting that under the assumptions (G1) and (G2), one can find some $C>0$ such that
$$g(s)\leq C\left(s^{\alpha-1}+s^{\beta-1}\right).$$
So the right-hand side  $\frac{1}{a+bd_n}\left[\frac{g(M_n U_n(y))}{M_{n}^{\beta-1}}-\frac{\lambda_n}{M_{n}^{\beta-2}}U_n(y)\right]$ is of $L^\infty(\R^N)$. Then applying a standard elliptic estimation and passing to a subsequence if necessary, $U_n\rightarrow U$ in $C_{loc}^{2}(\R^N)$,  which is a nontrivial and nonnegative bounded radial solution of
\beq\lab{eq:20210813-bxe2}
-\Delta U=\frac{B}{a} U^{\beta-1}~\hbox{in}~\R^N.
\eeq
Recalling the well known result  that the only nonnegative solution to $-\Delta u\geq u^p$ in $\R^N$ with $p<\frac{N}{(N-2)_+}$ is $0$, see \cite[Theorem 8.4]{QuittnerSouplet2007}. So for the case of $N=1,2$, we obtain that $U\equiv 0$, a contradiction to $U(0)=1$.
For the case of $N=3$, since $\beta-1<\frac{N+2}{N-2}$, by the remarkable result \cite{ChenLi1991}, the only nonnegative solution of \eqref{eq:20210813-bxe2} is also $0$, a contradiction to $U(0)=1$ again.
\ep

Define
\beq\lab{eq:20210810-we2}
\bar{u}_n(x):=\frac{1}{u_n(0)} u_n\left(\frac{x}{\sqrt{\lambda_n}}\right).
\eeq
A direct computation shows that $\bar{u}_n(0)=\|\bar{u}\|_\infty=1$ and
\begin{align}\lab{eq:20210810-we3}
&-\left(a +b d_{n}\right)\Delta \bar{u}_n(x) + \bar{u}_n(x)
=\frac{1}{\lambda_n u_n(0)}g(u_n(0)\bar{u}_n(x)).
\end{align}

\bl\lab{lemma:20210811-l1}
Under the assumptions of (G1)-(G5), let $(\lambda_c,u_c)$ be the solution given by Theorem \ref{thm:Main-Thoerem}. We have that
$$\liminf_{c\rightarrow +\infty}\frac{u_c(0)^{\alpha-2}}{\lambda_c}>0.$$
\el
\bp
We argue by contradiction and assume that there exists a sequence $c_n\rightarrow +\infty$ such that $u_n(0)^{\alpha-2}=o(1)(\lambda_n)$.
Consider $x=0$ in \eqref{eq:20210810-we3}, we have
\begin{align*}
1=\bar{u}_n(0)\leq &-\left(a+bd_{n}\right)\Delta \bar{u}_n(0)+\bar{u}_n(0)\\
=&\frac{1}{\lambda_{n}u_n(0)}g(u_n(0)\bar{u}_n(0))\\
\leq& C \frac{1}{\lambda_{n}u_n(0)}\left(u_n(0)^{\alpha-1}+u_n(0)^{\beta-1}\right)\\
\leq &C \frac{u_n(0)^{\alpha-2}}{\lambda_n}\\
=&o(1),
\end{align*}
a contradiction.
\ep

\bl\lab{cor:20210817-c1}
Under the assumptions of (G1)-(G5), let $(\lambda_c,u_c)$ be the solution given by Theorem \ref{thm:Main-Thoerem}. Then
$u_c(0)=\|u_c\|_\infty\rightarrow 0$ as $c\rightarrow +\infty$.
\el
\bp
Noting that $\lambda_c\rightarrow 0^+$ and $\|\nabla u_c\|_2\rightarrow 0$ as $c\rightarrow +\infty$, suppose now that $\displaystyle\liminf_{n\rightarrow +\infty}u_n(0)>0$. Since  $g(u_n)-\lambda_nu_n\in L^\infty(\R^N)$, applying the elliptic regularity and passing to a subsequence, we may assume that $u_n\rightarrow u$ in $C_{loc}^{2}(\R^N)$ with $u(0)=\max_{x\in \R^N}u(x)>0$ and
$u$  a non-negative bounded radial function which solves
$$-\Delta u=\frac{1}{a}g(u)\geq 0~\hbox{in}~\R^N.$$
Hence, by \cite[Theorem 7.2]{JeanZhangZhong2021}, we obtain that $u\equiv 0$, a contradiction to $u(0)>0$.
\ep

\bl\lab{lemma:20210811-l2}
Under the assumptions of (G1)-(G5),let $(\lambda_c,u_c)$ be the solution given by Theorem \ref{thm:Main-Thoerem}. We have that
$$\limsup_{c\rightarrow +\infty}\frac{u_c(0)^{\alpha-2}}{\lambda_c}<+\infty.$$
\el
\bp
If there exists a sequence $c_n\rightarrow +\infty$ such that $\frac{u_n(0)^{\alpha-2}}{\lambda_n}\rightarrow +\infty$. Then $\lambda_n=o\left(u_n(0)^{\alpha-2}\right)$.
Put
$$\tilde{u}_n(x):=\frac{1}{u_n(0)}u_n\left(\frac{1}{u_n(0)^{\frac{\alpha-2}{2}}}x\right),$$
then a direct computation shows that $\tilde{u}_n(0)=\|\tilde{u}_n\|_\infty=1$ and
\beq\lab{eq:20210817-xbue1}
-\left(a+bd_n\right)\Delta \tilde{u}_n=\frac{g(u_n(0)\tilde{u}_n)}{u_n(0)^{\alpha-1}}-\frac{\lambda_n}{u_n(0)^{\alpha-2}} \tilde{u}_n.
\eeq
Under the assumption (G2), by Lemma \ref{lemma:20210813-l1} and $\frac{\lambda_n}{u_n(0)^{\alpha-2}}=o(1)$, we see that the right hand side of \eqref{eq:20210817-xbue1} is of $L^\infty(\R^N)$. Hence, we can apply the standard elliptic estimation and suppose that $\tilde{u}_n\rightarrow \tilde{u}$ in $C_{loc}^{2}(\R^N)$, going up to a subsequence if necessary. Then by Lemma \ref{cor:20210817-c1} and the assumption (G4), combining with $d_n\rightarrow 0$, we see that $\tilde{u}$ is is a nonnegative radial bounded function satisfying
\beq\lab{eq:20210817-xbue2}
-\Delta \tilde{u}= \frac{A}{a}\tilde{u}^{\alpha-1}~\hbox{in}~\R^N.
\eeq
Then similar to the proof of Lemma \ref{lemma:20210813-l1}, we have that $\tilde{u}\equiv 0$, a contradiction to $\tilde{u}(0)=1$.
\ep

\bl\lab{lemma:20210811-xl1}
Under the assumptions (G1)-(G5), $\bar{u}_n(x)\rightarrow 0$ as $|x|\rightarrow +\infty$ uniformly in $n\in \N$.
\el
\bp
 We argue by contradiction and assume that there exist $\varepsilon>0$ and $r_n\rightarrow +\infty$ such that $\displaystyle\bar{u}_n(r_n)=\varepsilon <<\liminf_{c\rightarrow +\infty}\frac{u_c(0)^{\alpha-2}}{\lambda_c}$.
 By Lemma \ref{lemma:20210811-l1} and Lemma \ref{lemma:20210811-l2}, up to a subsequence, we may assume that
 $$\frac{u_n(0)^{\alpha-2}}{\lambda_n}\rightarrow C_0>>\varepsilon>0.$$
 Set $\tilde{u}_n(r)=\bar{u}_n(r+r_n)$, then by \eqref{eq:20210810-we3} we have that
 \beq\lab{eq:20210811-we1}
 -(a+bd_{n})\left({\tilde{u}_n}''+\frac{N-1}{r+r_n}{\tilde{u}_n}'\right)=-\tilde{u}_n+
 \frac{1}{\lambda_{n}u_n(0)}g(u_n(0)\tilde{u}_n(r)),\quad r>-r_n.
 \eeq
 By the standard elliptic estimation, for any compact subset  $\Omega\subset \R$, up to a subsequence, $\tilde{u}_n(r)\rightarrow \tilde{u}$ uniformly in $\Omega$. Then one can see that $\tilde{u}$ solves
 \beq\lab{eq:20210811-we2}
 -a \tilde{u}''=-\tilde{u}+C_0A\tilde{u}^{\alpha-1}, r\in \R,
 \eeq
 and $\tilde{u}(0)=\varepsilon$.
 Noting that $\tilde{u}_n$ is decreasing in $[0,+\infty)$, by $r_n\rightarrow +\infty$, we obtain that $\tilde{u}$ is decreasing in $\R$. Hence, $\tilde{u}$ has limit $\tilde{u}_+$ at $r=+\infty$ and also has limit $\tilde{u}_-$ at $r=-\infty$. In particular, $0\leq \tilde{u}_+<\tilde{u}_-\leq 1$. Hence, both $\tilde{u}_+$ and $\tilde{u}_-$ are solutions of $$-\tilde{u}+C_0A\tilde{u}^{\alpha-1}=0.$$
 Since $\varepsilon>0$ is very small, and $\tilde{u}_+\leq \tilde{u}(0)<\varepsilon$, one can see that $\tilde{u}_+=0$. On the other hand, $\tilde{u}_-\geq \tilde{u}(0)=\varepsilon>0$, we obtain that
 $\tilde{u}_-=(C_0A)^{-\frac{1}{\alpha-2}}$. Hence,
 $0\leq \tilde{u}(r)\leq (C_0A)^{-\frac{1}{\alpha-2}}$  in $\R$ and $0<\tilde{u}(r)\leq \varepsilon$ in $[0,+\infty)$.  Then we obtain that $-\tilde{u}''\leq 0$ in $\R$ and $-\tilde{u}''<0$ in $[0,+\infty)$. So $\tilde{u}$ is convex in $\R$ and strictly convex in $[0,+\infty)$.  A contradiction to $0\leq \tilde{u}(r)\leq (C_0A)^{-\frac{1}{\alpha-2}}$.
\ep

{\bf Proof of Theorem \ref{thm:20210814-th1}:}
For any sequence $c_n\rightarrow +\infty$,
we put
$$v_n(x):=\lambda_{n}^{\frac{1}{2-\alpha}} u_n\left(\frac{x}{\sqrt{\lambda_n}} \right)
=\frac{u_n(0)}{\lambda_{n}^{\frac{1}{\alpha-2}}} \bar{u}_n.$$
Under the assumptions (G1)-(G5), we remark that by Lemma \ref{lemma:20210811-l2},
$$\limsup_{n\rightarrow \infty}\sup_{x\in \R^N}v_n(x)=\limsup_{n\rightarrow \infty}\sup_{x\in \R^N}\lambda_{n}^{\frac{1}{2-\alpha}} u_n\left(\frac{x}{\sqrt{\lambda_n}} \right)=\limsup_{n\rightarrow \infty}\lambda_{n}^{\frac{1}{2-\alpha}} u_n(0)<+\infty.$$
Then combining with Lemma \ref{lemma:20210811-xl1}, we can see that $\displaystyle \lim_{|x|\rightarrow +\infty}v_n(x)=0$ uniformly in $n\in \N$. And $v_n$ solves the following equation
\beq\lab{eq:20210811-wbe1}
-(a+b\|\nabla u_n\|_2^2) \Delta v_n +v_n=\frac{g(\lambda_{n}^{\frac{1}{\alpha-2}} v_n)}{\lambda_{n}^{\frac{\alpha-1}{\alpha-2}}}.
\eeq
Hence, up to a subsequence, the standard elliptic estimate implies that $v_n(x)\rightarrow U(x)$ in $C_{loc}^{2}(\R^N)$.
Recalling that $\|\nabla u_n\|_2^2\rightarrow 0, \lambda_n\rightarrow 0$, by (G4), \eqref{eq:20210811-wbe1} implies that   $U$ is nontrivial nonnegative function satisfying
$$\begin{cases}
-a\Delta U+U=A U^{\alpha-1}~\hbox{in}~\R^N,\\
\lim_{|x|\rightarrow +\infty}U(x)=0.
\end{cases}$$
Then $U$ is uniquely determined by the well-known result of Kwong \cite{Kwong1989}, which is a positive radial decreasing function.
\hfill$\Box$

{\bf Proof of Theorem \ref{thm:20210818-th1}:} Recalling \eqref{eq:20210811-wbe1}, we have
\beq\lab{eq:20210818-ze1}
-\Delta v_n +\frac{1}{a+b\|\nabla u_n\|_2^2} \left[1 -\frac{g\left(\lambda_{n}^{\frac{1}{\alpha-2}}v_n\right)}{\lambda_{n}^{\frac{\alpha-1}{\alpha-2}}v_n}\right]v_n=0.
\eeq
Noting that $\|v_n\|_\infty=1$ and $\lambda_n\rightarrow 0^+$, under the assumption (G4), we have
\beq\lab{eq:20210818-ze2}
g\left(\lambda_{n}^{\frac{1}{\alpha-2}}v_n(x)\right)=\left(A+o(1)\right) \lambda_{n}^{\frac{\alpha-1}{\alpha-2}}v_{n}(x)^{\alpha-1}~\hbox{as}~n\rightarrow +\infty.
\eeq
Then for $n$ large, combing with $\|\nabla u_n\|_2^2=d_n\rightarrow 0$, we can write \eqref{eq:20210818-ze1} as
\beq\lab{eq:20210818-ze3}
-\Delta v_n +\left(\frac{1}{a}+o(1)\right)\left[1-\left(A+o(1)\right)v_{n}^{\alpha-2}\right]v_n=0.
\eeq
Recalling that $\displaystyle \lim_{|x|\rightarrow +\infty}v_n(x)=0$ uniformly in $n\in \N$, we can find some $R$ large enough such that
\beq\lab{eq:20210818-ze4}
\left(\frac{1}{a}+o(1)\right)\left[1-\left(A+o(1)\right)v_{n}^{\alpha-2}\right]>\frac{1}{2a}, \forall |x|\geq R~\hbox{and}~n\geq N_0.
\eeq
Hence,
\beq\lab{eq:20210818-ze5}
-\Delta v_n(x)+\frac{1}{2a}v_n(x)\leq 0, \forall |x|\geq R, \forall n\geq N_0.
\eeq
On the other hand, by $v_n(x)\rightarrow U(x)$ in $C_{loc}^{2}(\R^N)$,  one can find some $C_1,C_2>0$ such that
\beq\lab{eq:20210818-ze6}
v_n(x)\leq C_1 e^{-C_2|x|}, \forall |x|\geq R, \forall n\geq N_0.
\eeq
Hence, $v_n\rightarrow U$ in $L^2(\R^N)$.
A direct computation shows that
$$\|v_n\|_2^2=\lambda_{n}^{\frac{N\alpha-2N-4}{2(\alpha-2)}}c_n~\hbox{and}~\|\nabla v_n\|_2^2=\lambda_{n}^{\frac{N\alpha-2N-4}{2(\alpha-2)}} \frac{d_n}{\lambda_n}.$$
Recalling Remark \ref{remark:20210812-r1} and noting that $d_n\rightarrow 0$, we obtain that $\|v_n\|_2^2$ and $\|\nabla v_n\|_2^2$ are comparable. By $v_n\rightarrow U$ in $L^2(\R^N)$, there exists some $C_3,C_4>0$ such that
$$C_3\leq \|v_n\|_2^2\leq C_4, \forall n\in \N.$$
And thus we can also find some $C_5,C_6>0$ such that
$$C_5\leq \|\nabla v_n\|_2^2\leq C_6.$$
So $\{v_n\}$ is a bounded sequence of $H^1(\R^N)$. Recalling $v_n\rightarrow U$ in $L^2(\R^N)$ again, up to a subsequence, we can have that
$$v_n\rightharpoonup U~\hbox{in}~H^1(\R^N).$$
Furthermore,  by Gagliardo-Nierenberg inequality, one can see that
$$
v_n\rightarrow U~\hbox{in}~L^p(\R^N), \forall p\in [2, 2^*).
$$
Then under the assumption (G2), by \eqref{eq:20210818-ze3}, we have
$$\|\nabla v_n\|_2^2\rightarrow \|\nabla U\|_2^2.$$
So we also have that $v_n\rightarrow U$ in $D_{0}^{1,2}(\R^N)$.
Hence, $v_n\rightarrow U$ in $H^1(\R^N)$.
\hfill$\Box$

\subsection{The case of $c\rightarrow 0^+$}
Let $c_n\rightarrow 0^+$, by Corollary \ref{cor:20210810-c1}, we have that  $\lambda_{c_n}\rightarrow +\infty$. Write $(\lambda_{c_n}, u_{c_n})$ simply by $(\lambda_n, u_n)$.
Also put $d_n:=\|\nabla u_n\|_2^2$.
And by \eqref{eq:20210810-e8}, we have that $d_n\rightarrow +\infty$ as $n\rightarrow +\infty$.

\bl\lab{lemma:20210814-l1}
Under the assumptions (G1)-(G5), let $(\lambda_c,u_c)$ be the solution given by Theorem \ref{thm:Main-Thoerem}. Then
$$\liminf_{c\rightarrow 0^+} u_c(0)=+\infty$$
and
$$\liminf_{c\rightarrow 0^+} \frac{u_c(0)^{\beta-2}}{\lambda_c}>0.$$
\el
\bp
Define
$$v_n(x):=\frac{1}{u_n(0)}u_n\left(\frac{\sqrt{d_n}}{\sqrt{\lambda_n}}x\right).$$
Then a direct computation shows that $\|v_n\|_\infty=v_n(0)=1$ and
\beq\lab{eq:20210814-xe1}
-\left(ad_{n}^{-1}+b\right)\Delta v_n+v_n=\frac{1}{\lambda_n u_n(0)} g(u_n(0) v_n).
\eeq
By (G2) and taking $x=0$, there exists some $C>0$ such that
\begin{align*}
1\leq & -\left(ad_{n}^{-1}+b\right)\Delta v_n (0)+v_n(0)=\frac{1}{\lambda_n u_n(0)} g(u_n(0) v_n(0))\\
\leq& C\frac{1}{\lambda_n u_n(0)}\left(u_n(0)^{\alpha-1}+u_n(0)^{\beta-1}\right)\\
=&C\frac{u_n(0)^{\alpha-2}+u_n(0)^{\beta-2}}{\lambda_n}.
\end{align*}
Since $\lambda_n\rightarrow +\infty$, we obtain that $u_n(0)\rightarrow +\infty$ since $\alpha-2>0,\beta-2>0$. In particular, by $\alpha\leq \beta$, we obtain that
$$\liminf_{n\rightarrow +\infty}\frac{u_n(0)^{\beta-2}}{\lambda_n}\geq \frac{1}{2C}>0.$$
By the arbitrary of $c_n$, we obtain $\liminf_{c\rightarrow 0^+} u_c(0)=+\infty$ and
$$\liminf_{c\rightarrow 0^+} \frac{u_c(0)^{\beta-2}}{\lambda_c}>0.$$
\ep

\bl\lab{lemma:20210814-l2}
Under the assumptions (G1)-(G5), let $(\lambda_c,u_c)$ be the solution given by Theorem \ref{thm:Main-Thoerem}. Then
$$\limsup_{c\rightarrow 0^+} \frac{u_c(0)^{\beta-2}}{\lambda_c}<+\infty.$$
\el
\bp
We argue by contradiction and suppose there exists a sequence $c_n\rightarrow 0^+$ such that
$$\frac{u_n(0)^{\beta-2}}{\lambda_n}\rightarrow +\infty.$$
That is, $\lambda_n=o(u_n(0)^{\beta-2})$.
Now, we put
\beq\lab{eq:20210814-xe2}
\bar{v}_n(x):=\frac{1}{u_n(0)}u_n\left(\sqrt{\frac{d_n}{u_n(0)^{\beta-2}}}x\right).
\eeq
Then $\bar{v}_n(0)=\|\bar{v}\|_\infty=1$ and
\beq\lab{eq:20210814-xe3}
-\left(ad_n^{-1}+b\right)\Delta \bar{v}_n+\frac{\lambda_n}{u_n(0)^{\beta-2}}\bar{v}_n=\frac{1}{u_n(0)^{\beta-1}} g(u_n(0)\bar{v}_n).
\eeq
Noting that
\begin{align*}
\frac{1}{u_n(0)^{\beta-1}} g(u_n(0)\bar{v}_n)\leq &C \frac{1}{u_n(0)^{\beta-1}} \left(u_n(0)^{\alpha-1}\bar{v}_{n}^{\alpha-1}+u_n(0)^{\beta-1}\bar{v}_{n}^{\beta-1}\right)\\
=&C \left(u_n(0)^{\alpha-\beta} \bar{v}_{n}^{\alpha-1}+\bar{v}_{n}^{\beta-1}\right).
\end{align*}
Recalling Lemma \ref{lemma:20210814-l1}, $u_n(0)\rightarrow +\infty$. By $\alpha\leq \beta$, we obtain that $$\frac{1}{u_n(0)^{\beta-1}} g(u_n(0)\bar{v}_n)\leq C \left(\bar{v}_{n}^{\alpha-1}+\bar{v}_{n}^{\beta-1}\right),$$
which is of $L^\infty(\R^N)$. Hence, we can apply the standard elliptic estimation and assume that $\bar{v}_n\rightarrow \bar{v}$ in $C_{loc}^{2}(\R^N)$. So $\bar{v}$ is a nonnegative radial bounded function and (G5) implies that
\beq\lab{eq:20210814-xe4}
-\Delta \bar{v}=\frac{B}{b}\bar{v}^{\beta-1}~\hbox{in}~\R^N.
\eeq
Then $\bar{v}\equiv 0$, a contradiction to $\bar{v}(0)=1$, a similar argument we refer to the proof of Lemma \ref{lemma:20210813-l1}.
\ep

{\bf Proof of Theorem \ref{thm:20210814-th2}:}
For any sequence $c_n\rightarrow 0^+$, a direct computation shows that $v_n$ satisfies
\beq\lab{eq:20210814-xe5}
-\left(a d_{n}^{-1}+b\right)\Delta v_n+v_n=\frac{g (\lambda_{n}^{\frac{1}{\beta-2}} v_n)}{\lambda_{n}^{\frac{\beta-1}{\beta-2}}}.
\eeq
By Lemma \ref{lemma:20210814-l2}, there exists some $M>0$ such that
$$v_n(0)=\|v_n\|_\infty\leq M, \forall n\in \N.$$
So combining with the assumption (G2), we have
\begin{align*}
\frac{g (\lambda_{n}^{\frac{1}{\beta-2}} v_n)}{\lambda_{n}^{\frac{\beta-1}{\beta-2}}}
\leq& \frac{C}{\lambda_{n}^{\frac{\beta-1}{\beta-2}}}\left(\lambda_{n}^{\frac{\alpha-1}{\beta-2}} v_{n}^{\alpha-1}+\lambda_{n}^{\frac{\beta-1}{\beta-2}}v_{n}^{\beta-1}\right)\\
=&C \left(\lambda_{n}^{\frac{\alpha-\beta}{\beta-2}} v_{n}^{\alpha-1}+v_{n}^{\beta-1}\right).
\end{align*}
So by $\alpha\leq \beta$ and $\lambda_n\rightarrow +\infty$, we see that $\frac{g (\lambda_{n}^{\frac{1}{\beta-2}} v_n)}{\lambda_{n}^{\frac{\beta-1}{\beta-2}}}\in L^\infty(\R^N)$.
 So applying a standard elliptic estimation, we can prove that $v_n\rightarrow v$ in $C_{loc}^{2}(\R^N)$ up to a subsequence, where $v$ is a bounded nonnegative radial decreasing function satisfying
\beq\lab{eq:20210814-xe6}
-b\Delta v+v=B v^{\beta-1}~\hbox{in}~\R^N.
\eeq
We claim that $\displaystyle \lim_{|x|\rightarrow +\infty}v_n(x)=0$ uniformly in $n\in \N$.
If not, up to a subsequence, we suppose there exist some $\varepsilon>0$ small enough and sequence $r_n\rightarrow +\infty$ such that $v_n(r_n)=\varepsilon>0$. We put $\tilde{v}_n(r):=v_n(r+r_n)$. Then  by \eqref{eq:20210814-xe5}, we have that
\beq\lab{eq:20210814-xe7}
-\left(a d_{n}^{-1}+b\right) \left(\tilde{v}''_n+\frac{N-1}{r+r_n}\tilde{v}'_n\right)+\tilde{v}_n=\frac{g (\lambda_{n}^{\frac{1}{\beta-2}} v_n)}{\lambda_{n}^{\frac{\beta-1}{\beta-2}}}, r\geq -r_n.
\eeq
By the standard elliptic estimation, for any compact subset  $\Omega\subset \R$, up to a subsequence, $\tilde{v}_n(r)\rightarrow \tilde{v}$ uniformly in $\Omega$. Then one can see that $\tilde{v}$ solves
 \beq\lab{eq:20210814-xe8}
 -b \tilde{v}''=-\tilde{v}+B\tilde{v}^{\beta-1}, r\in \R,
 \eeq
 and $\tilde{v}(0)=\varepsilon$.
 Noting that $\tilde{v}_n$ is decreasing in $[0,+\infty)$, by $r_n\rightarrow +\infty$, we obtain that $\tilde{v}$ is decreasing in $\R$. Hence, $\tilde{v}$ has limit $\tilde{v}_+$ at $r=+\infty$ and also has limit $\tilde{v}_-$ at $r=-\infty$. In particular, $0\leq \tilde{v}_+<\tilde{v}_-\leq \limsup_{c\rightarrow 0^+}\frac{u_{c}(0)}{\lambda_{c}^{\frac{1}{\beta-2}}}<+\infty$. Hence, both $\tilde{v}_+$ and $\tilde{v}_-$ are solutions of $$-\tilde{v}+B\tilde{v}^{\beta-1}=0.$$
 Since $\varepsilon>0$ is very small, and $\tilde{v}_+\leq \tilde{v}(0)=\varepsilon$, one can see that $\tilde{v}_+=0$. On the other hand, $\tilde{v}_-\geq \tilde{v}(0)=\varepsilon>0$, we obtain that
 $\tilde{v}_-=B^{-\frac{1}{\beta-2}}$. Hence,
 $0\leq \tilde{v}(r)\leq B^{-\frac{1}{\beta-2}}$  in $\R$ and $0<\tilde{v}(r)\leq \varepsilon$ in $[0,+\infty)$.  Then we obtain that $-\tilde{v}''\leq 0$ in $\R$ and $-\tilde{v}''<0$ in $[0,+\infty)$. So $\tilde{v}$ is convex in $\R$ and strictly convex in $[0,+\infty)$.  A contradiction to $0\leq \tilde{v}(r)\leq B^{-\frac{1}{\beta-2}}$.
 Hence, we prove that $\displaystyle \lim_{|x|\rightarrow +\infty}v_n(x)=0$ uniformly in $n\in \N$ and thus
 \beq\lab{eq:20210814-xe9}
 \lim_{|x|\rightarrow +\infty}v(x)=0.
 \eeq
 That is, $v$ is a nontrivial nonnegative function satisfying \eqref{eq:20210814-xe6} and \eqref{eq:20210814-xe9}, then by the classification result of Kwong \cite{Kwong1989}, $v$ is positive radial decreasing and uniquely determined.
\hfill$\Box$

{\bf Proof of Theorem \ref{thm:20210818-th2}:}
Let $c_n\rightarrow 0^+$, recalling \eqref{eq:20210814-xe5},
\beq\lab{eq:20210818-we1}
-\Delta v_n+\frac{1}{a d_{n}^{-1}+b}\left[1-\frac{g (\lambda_{n}^{\frac{1}{\beta-2}} v_n)}{\lambda_{n}^{\frac{\beta-1}{\beta-2}}v_n}\right]v_n=0.
\eeq
By the proof of Theorem \ref{thm:20210814-th2}, we also have that $\displaystyle \lim_{|x|\rightarrow +\infty}v_n(x)=0$ uniformly in $n\in \N$.
So by (G2) and $\lambda_n\rightarrow +\infty$ as $n\rightarrow +\infty$,
\begin{align*}
&\frac{g (\lambda_{n}^{\frac{1}{\beta-2}} v_n)}{\lambda_{n}^{\frac{\beta-1}{\beta-2}}v_n}
\leq C \frac{\lambda_{n}^{\frac{\alpha-1}{\beta-2}}v_{n}^{\alpha-1}+\lambda_{n}^{\frac{\beta-1}{\beta-2}}v_{n}^{\beta-1}}
{\lambda_{n}^{\frac{\beta-1}{\beta-2}}v_n}\\
\leq&C\left[\lambda_{n}^{\frac{\alpha-\beta}{\beta-2}}v_{n}^{\alpha-1}+v_{n}^{\beta-1}\right]
\leq C\left[v_{n}^{\alpha-1}+v_{n}^{\beta-1}\right]\\
\rightarrow &0~\hbox{as}~|x|\rightarrow +\infty~\hbox{uniformly in $n\in \N$}.
\end{align*}
Hence, we can find some $R>0$ large enough and some $N_0\in \N$ such that
\beq\lab{eq:20210818-we2}
\frac{1}{a d_{n}^{-1}+b}\left[1-\frac{g (\lambda_{n}^{\frac{1}{\beta-2}} v_n)}{\lambda_{n}^{\frac{\beta-1}{\beta-2}}v_n}\right]\geq \frac{1}{2b}, \forall |x|\geq R~\hbox{and}~\forall ~n\geq N_0.
\eeq
Then applying a similar argument as the proof of Theorem \ref{thm:20210818-th1}, we can prove that $v_n\rightarrow V$ in $L^2(\R^N)$.
Noting that $d_n\rightarrow +\infty$, by Remark \ref{remark:20210812-r1}, we see that $\lambda_n c_n$ and $d_n^2$ are comparable.
Since
$$\|v_n\|_2^2=\lambda_{n}^{\frac{N}{2}-\frac{2}{\beta-2}} d_{n}^{-\frac{N}{2}} c_n, \|\nabla v_n\|_2^2=\lambda_{n}^{\frac{N}{2}-\frac{2}{\beta-2}} d_{n}^{-\frac{N}{2}} \frac{d_n^2}{\lambda_n},$$
we get that $\|v_n\|_2^2$ and $\|\nabla v_n\|_2^2$ are comparable.
Then similar to the proof of Theorem \ref{thm:20210818-th1}, we can finally that obtain $v_n\rightarrow V$ in $H^1(\R^N)$.
\hfill$\Box$

\end{document}